\documentclass[number,citesort,dvips]{arxbj}
\usepackage{upgreek}
\usepackage{graphicx}
\usepackage{breakurl}

\aid{0}
\volume{17}
\issue{2}
\pubyear{2011}
\firstpage{736}
\lastpage{748}
\doi{10.3150/10-BEJ296}

\makeatletter

\newtheorem{lemma}{Lemma}

\newtheorem{theorem}{Theorem}

\newproclaim{assu}{Assumptions}

\makeatother

\begin{document}
\begin{frontmatter}

\title{Semi-parametric regression: Efficiency gains from
modeling the nonparametric part}

\runtitle{Semi-parametric regression}

\begin{aug}
\author[1]{\fnms{Kyusang} \snm{Yu}\thanksref{1}\ead[label=e1]{kyusangu@konkuk.ac.kr}},
\author[2]{\fnms{Enno} \snm{Mammen}\thanksref{2}\ead[label=e2]{emammen@rumms.uni-mannheim.de}}
\and
\author[3]{\fnms{Byeong U.} \snm{Park}\corref{}\thanksref{3}\ead[label=e3]{bupark@stats.snu.ac.kr}}

\runauthor{K. Yu, E. Mammen and B.U. Park}

\address[1]{Konkuk University, Seoul, Korea. \printead{e1}}
\address[2]{University of Mannheim, Mannheimun, Germany. \printead{e2}}
\address[3]{Seoul National University, Seoul, Korea. \printead{e3}}
\end{aug}

\received{\smonth{4} \syear{2010}}

%
\begin{abstract}
It is widely admitted that structured nonparametric modeling that
circumvents the curse of dimensionality is important in nonparametric
estimation. In this paper we show that the same holds for
semi-parametric estimation. We argue that estimation of the parametric
component of a semi-parametric model can be improved essentially when
more structure is put into the nonparametric part of the model. We
illustrate this for the partially linear model, and investigate
efficiency gains when the nonparametric part of the model has an
additive structure. We present the semi-parametric Fisher information
bound for estimating the parametric part of the partially linear
additive model and provide semi-parametric efficient estimators for
which we use a smooth backfitting technique to deal with the additive
nonparametric part. We also present the finite sample performances of
the proposed estimators and analyze Boston housing data as an
illustration.
\end{abstract}

%

\end{frontmatter}

\section{Introduction}

Structured nonparametric models such as additive models are known to
circumvent the curse of dimensionality and allow reliable estimation
when a full nonparametric model does not work. In the present paper we
show that a similar assertion applies for semi-parametric models:
structural modeling of the nonparametric part can lead to accurate
estimation of the parametric part even in situations where otherwise
only very poor, unreliable or unstable estimates would be available. We
show this by comparing the partially linear and the partially linear
additive model. In particular, we demonstrate that using an additive
model for the nonparametric part in the partially linear model can lead
to drastic gains of efficiency in the estimation of the parametric
components. This holds if the dimension of the nonparametric covariates
is high, or the parametric covariates can be approximated by
non-additive transformations of the nonparametric covariates. In the
extreme of the latter case, if the approximation is exact, then
estimation of the parametric part in the partially linear model breaks
down. If the approximation is very crude, one sees large efficiency
gains by using additive models for the nonparametric part.\looseness=-1

Suppose we observe the i.i.d. copies
$(Y^1,\mathbf{X}^{1},\mathbf{Z}^{1}),\ldots,(Y^n,\mathbf
{X}^{n},\mathbf{Z}^{n})$ of a random vector
$(Y,\mathbf{X},\mathbf{Z})$, where $\mathbf{X}=(X_1,\ldots
,X_p)^\top\in\mathbb{R}^p$ and
$\mathbf{Z}=(Z_1,\ldots,Z_d)^\top\in\mathbb{R}^d$. The partially linear
model assumes
\begin{equation}\label{plmodel}
Y = m_0+\mathbf{X}^\top{\bolds{\beta}}+ m(Z_1,\ldots,Z_d)+\epsilon,
\end{equation}
where $\bolds{\beta}$ is an unknown $p$-vector and $m$ is an unknown
$d$-variate function. The partially linear additive model puts an
additive structure to the nonparametric function $m$:
\begin{equation}\label{pladdmodel}
Y = m_0+\mathbf{X}^\top\bolds{\beta}+ m_1(Z_1)+\cdots
+m_d(Z_d)+\epsilon.
\end{equation}
These models exclude the interesting case where $\mathbf{X}$ or $\mathbf{Z}$
includes some endogeneous variables of~$Y$, but they simplify our
discussion on semi-parametric efficiency. We believe that our results
can be extended to the corresponding semi-parametric models with time
series data by following, for example, the arguments in \cite{ks1997}.

For identifiability of the additive component functions $m_j$, we put
the constraints $E m_j(Z_j) = 0, ~1 \le j \le d$. We assume that
$(\mathbf{X},\mathbf{Z})$ has a joint density $q$ with respect to
$\nu=\nu_1
\times
\nu_2$, where $\nu_1$ is a $\sigma$-finite measure and $\nu_2$ is the
Lebesgue measure on each support of $\mathbf{X}$ and $\mathbf{Z}$,
and that the
marginal density of $\mathbf{Z}$ (with respect to $\nu_1$), denoted by
$q_{\mathbf{Z}}$, has compact support, say $[0,1]^d$. The model
(\ref{pladdmodel}) enjoys the advantages of both the partially linear
model (\ref{plmodel}) and the nonparametric additive model to the fully
nonparametric model. It accommodates discrete covariates since we only
require that $\nu_1$ is a $\sigma$-finite measure, and also interaction
effects between covariates by putting them into the parametric part. By
the additive structure in the nonparametric part it avoids the curse of
dimensionality, but retains the flexibility of the model. It also
renders easy interpretation of the individual role of each covariate.

We discuss semi-parametric efficient estimation of the parameter
$\bolds{\beta}$ in the model (\ref{pladdmodel}). We present the
semi-parametric Fisher information bound and provide an estimator that
achieves the efficiency bound. Semi-parametric efficient estimation
when $d=1$ has been studied by Bhattacharya and Zhao \cite{bz1997}, Cuzick
\cite{c1992} and Schick \cite{s1993}. Their works can be easily extended to the
model (\ref{plmodel}) for $d>1$. Comparing the Fisher information
bounds for the models (\ref{plmodel}) and (\ref{pladdmodel}), we find
that the information bound under the model (\ref{pladdmodel}) is
smaller than the bound under the model~(\ref{plmodel}). In our
semi-parametric model (\ref{pladdmodel}), we do not specify the
distribution of the error term $\epsilon$ or the distribution $q$ of
the covariates. We show that one can do as well without knowing those
distributions.

There have been a few works on the model (\ref{pladdmodel}). Opsomer
and Ruppert \cite{or1999} obtained a $\sqrt n$-consistent estimator of
$\bolds{\beta}$ by a backfitting method with undersmoothing. Recently Liang
\textit{et al.} \cite{ltrah2008} and Carroll \textit{et al.} \cite{cmmy2009} studied the model
with measurement error and repeated measurements, respectively. But
they did not discuss semiparametric efficiency. The model
(\ref{plmodel}) has been studied more often; see \cite{s1988}, among
others. Most studies, however, are rather focused on the cases where
there is only a single-dimensional (or at most low-dimensional)
nonparametric function $m$. This is because high-dimension costs
higher-order smoothness in theory and poor small sample performances in
practice.

\section{Semi-parametric efficiency}\label{sec:eff}

To avoid unnecessary complexity, we assume $m_0=0$. We also assume that
$\epsilon$ is independent with $(\mathbf{X},\mathbf{Z})$, and that
$g$, the density
of $\epsilon$, is symmetric and is absolutely continuous with respect
to the Lebesgue measure, having a derivative $g'$ and finite Fisher
information $\int(g')^2/g < \infty$. Below, we give a heuristic
argument for deriving the semi-parametric efficiency and present a
rigorous statement in a theorem.

Suppose that $g$ is known and $p=1$. We write $m(\mathbf{z}) = m_1(z_1) +
\cdots+ m_d(z_d)$ and adopt the convention $m_j(\mathbf{z})=m_j(z_j)$. The
logarithm of the joint density of $(Y,\mathbf{X},\mathbf{Z})$ as a
function of the
parameters is given by $\ell(\beta, m; (y,x,\mathbf{z})) = \log g
(y-x\beta-
m(\mathbf{z}))$, neglecting those terms that do not depend on $(\beta, m)$,
and the log-likelihood of $(\beta,m)$ by $\sum_{i=1}^n
\ell(\beta,m;(Y^i,X^i,\mathbf{Z}^i))$. Let $\mathcal{H}$ denote the
space of all
additive functions $m$ such that $m(\mathbf{z})=m_1(z_1)+\cdots
+m_d(z_d)$, $E
m_j(Z_j)=0$ and $Em(\mathbf{Z})^2 < \infty$.

Calculation of the Fisher information in a semi-parametric model is
made locally: fix a value $(\beta^0, m^0)$ of the parameter $(\beta,
m)$ and think of all `regular' parametric submodels $\{(\beta,
m_\beta)\dvtx \beta\in\mathbb{R}\}$ passing through $(\beta^0, m^0)$,
where $m_{\beta^0}=m^0$ and the mapping $\beta\mapsto m_{\beta}$ is
Fr\'echet differentiable as a function from $\mathbb{R}$ to $\mathcal{H}$.
Define $\varphi=g'/g$. Then, each finite-dimensional submodel
$\{(\beta, m_\beta)\dvtx\beta\in\mathbb{R}\}$ has the score function
\begin{eqnarray*}
\mathrm{d}\ell(\beta,m_{\beta})/\mathrm{d}\beta|_{\beta=\beta^0} &=&
\partial\ell(\beta,m^0)/\partial\beta|_{\beta=\beta^0}
+\partial\ell(\beta^0,m)/\partial m |_{m=m^0}(\delta)\\
&=&\varphi(\epsilon) X + \varphi(\epsilon) \delta(\mathbf{Z}),
\end{eqnarray*}
where $\displaystyle\delta={\partial m_{\beta}}/
{\partial\beta}|_{\beta=\beta^0} \in\mathcal{H}$ is the
tangent of the
mapping $\beta\mapsto m_{\beta}$ at $\beta^0$, and ${\partial
\ell}/{\partial m}$ denotes the Fr\'{e}chet derivative of $\ell$ with
respect to $m$. This gives the Fisher information for estimating
$\beta$ in each submodel as $\mathcal{I}(\delta) \equiv E [\varphi
(\epsilon
)X +
\varphi(\epsilon) \delta(\mathbf{Z}) ]^2$.

The Fisher information at $(\beta^0, m^0) \in\mathbb{R} \times
\mathcal{H}$ in
the \textit{full} semi-parametric model typically equals to the Fisher
information at $(\beta^0, m^0) \in\mathbb{R} \times\mathcal{H}$ in
the most
difficult parametric submodel that gives minimal $\mathcal{I}(\delta)$.
Theorem \ref{thm:info} below demonstrates that this is the case with our problem.
The least favorable direction $\delta^*$ that minimizes $\mathcal
{I}(\delta)$
over $\delta\in\mathcal{H}$ is the solution of the following integral
equation: for all $\delta\in\mathcal{H},$
\begin{eqnarray*}
0 &=& E [\varphi(\epsilon)X + \varphi(\epsilon) \delta^*(\mathbf{Z})
]\varphi(\epsilon) \delta(\mathbf{Z})\\
&=& I_g \cdot E \bigl[\bigl(E(X|\mathbf{Z}) + \delta^*(\mathbf{Z})\bigr)
\delta(\mathbf{Z})\bigr],
\end{eqnarray*}
where $I_g = \int(g')^2/g$. This shows that $\delta^* = -
\Pi(E(X|\mathbf{Z}=\cdot) | \mathcal{H})$, where $\Pi(\cdot|
\mathcal{H})$ denotes the
projection operator onto $\mathcal{H}$, and that the `curve' $m_\beta^*$
corresponding to the least favorable submodel equals $m_\beta^* =
(\beta^0 - \beta)\Pi(E(X|\mathbf{Z}=\cdot) | \mathcal{H}) +m^0$.
The Fisher
information for the least favorable submodel is thus given by
$\mathcal{I}(\delta^*) = I_g \cdot E [X-\Pi(E(X|\mathbf
{Z})|\mathcal{H})]^2$, where, with a
slight abuse of notation, we write $\Pi(E(X|\mathbf{Z}=\cdot
)|\mathcal{H})(\mathbf{Z}
)=\Pi
(E(X|\mathbf{Z})|\mathcal{H})$.

The above arguments can be generalized to the case where $p>1$. Writing
$\eta_j = \Pi(E(X_j|\mathbf{Z}=\cdot)|\mathcal{H})$ and $\bolds
{\eta}= (\eta_1, \ldots,
\eta_p)^\top$, the least favorable direction equals $\bolds{\delta
}^* = -
\bolds{\eta}$ so that the Fisher information matrix for the least favorable
submodel equals $\mathcal{I}(\bolds{\delta}^*) = I_g \cdot
E[\mathbf{X}-\bolds{\eta}(\mathbf{Z})][\mathbf{X}-\bolds{\eta
}(\mathbf{Z})]^\top$. In the following
theorem we
show that the Fisher information $\mathcal{I}(\bolds{\delta}^*)$
given above is indeed
the semi-parametric information bound, as defined in \cite{bkrw1993}, in our original semi-parametric model where the error
density $g$ and the density $q$ of the covariate $(\mathbf{X},\mathbf
{Z})$ are not
specified. To state the theorem, let $\mathcal{G}$ denote the set of all
symmetric and absolutely continuous (with respect to the Lebesgue
measure) functions $g$ such that $I_g <\infty$. Let $\mathcal{Q}$ be an
arbitrary class of density functions $q$. For the spaces of $m$, we
consider Hilbert spaces defined by
\[
\mathcal{H}(q) = \Biggl\{m \in L_2(q)\dvtx m(\mathbf{z})= \sum_{j=1}^d
m_j(z_j) \mbox
{ and
} E m_j(Z_j) =0 \mbox{ for all } 1 \le j \le d \Biggr\},
\]
where $L_2(q)$ denotes the space of functions $m\dvtx \mathbb{R}^d
\rightarrow\mathbb{R}$ such that $E_q m(\mathbf{Z})^2 <\infty$ and $E_q$
means the expectation under the density $q$. The semi-parametric model
(\ref{pladdmodel}) under study is then expressed as $\mathcal{P}= \{
p(\cdot;
\bolds{\beta}, m, g, q) \dvtx \bolds{\beta}\in\mathbb{R}^p, m \in
\mathcal{H}(q), g \in\mathcal{G}
, q
\in\mathcal{Q}\}$. Let $(\bolds{\beta}^0, m^0, g_0, q_0)$ be a
fixed point where we
are calculating the semi-parametric Fisher information. Denote by $P_0$
the distribution corresponding to $(\bolds{\beta}^0, m^0, g_0, q_0)$,
and by
$I(P_0 | \bolds{\beta}, \mathcal{P})$ the semi-parametric Fisher
information at $P_0$
for estimating $\bolds{\beta}$ under the model $\mathcal{P}$. In the
theorem below,
the `efficient score' $\ell^*$ for estimating $\bolds{\beta}$ is the
score for
$\bolds{\beta}$ at $\bolds{\beta}^0$ in the least favorable
parametric submodel that
is indexed only by $\bolds{\beta}$ and passes through $P_0$. Let
$E_0$ denote
the expectation under $P_0$.

\begin{theorem}\label{thm:info}
The efficient score at $P_0$ for estimating $\bolds{\beta}$ is given by
\begin{eqnarray*}
&&\ell^*(\mathbf{x}, \mathbf{z}, y; P_0 | \bolds{\beta}, \mathcal
{P})
\\
&&\quad
=-[\mathbf{x}- \bolds{\eta}
(\mathbf{z})]
\frac{{g_0}'}{g_0}\bigl(y-\mathbf{x}^\top\bolds{\beta}^0 - m^0(\mathbf{z})\bigr),
\end{eqnarray*}
where $\bolds{\eta}= (\Pi[E_0(X_j |\mathbf{Z}=\cdot) | \mathcal
{H}(q_0)])_{j=1}^p$. The
information bound at $P_0$ for estimating $\bolds{\beta}$ equals
$I(P_0 |
\bolds{\beta}, \mathcal{P}) = I_{g_0} \cdot E_0 [\mathbf{X}-\bolds
{\eta}(\mathbf{Z})]
[\mathbf{X}-\bolds{\eta}(\mathbf{Z})]^\top$.
\end{theorem}

A proof of Theorem~\ref{thm:info} can be found in an extended version
of this paper that can be downloaded from \url{http://stat.snu.ac.kr/theostat/papers/BEJ296\_ExtendedVersion.pdf}.

Let $\mathcal{P}_{\mathrm{PL}} \supset\mathcal{P}$ denote the
semi-parametric model
(\ref{plmodel}). One can show $I(P_0 | \bolds{\beta}, \mathcal
{P}_{\mathrm{PL}})=I_{g_0}
\cdot E_0 [\mathbf{X}-E_0(\mathbf{X}|\mathbf{Z})][\mathbf
{X}-E_0(\mathbf{X}|\mathbf{Z})]^\top$ using the
arguments to derive $I(P_0 | \bolds{\beta}, \mathcal{P})$. Note that
$I(P_0 | \bolds{\beta},
\mathcal{P}) \ge I(P_0 | \bolds{\beta}, \mathcal{P}_{\mathrm{PL}})$ by
the property of conditional
expectation, and that the equality $I(P_0 | \bolds{\beta}, \mathcal
{P}) = I(P_0 |
\bolds{\beta}, \mathcal{P}_{\mathrm{PL}})$ holds if $E_0(X_j|\mathbf
{Z}=\mathbf{z})$ are additive
for all
$1 \le j \le d$. According to the theory of semi-parametric efficiency,
the minimal asymptotic variance that any regular estimator of $\bolds
{\beta}$
can achieve equals the inverse of the Fisher information matrix. The
inequality $I(P_0 | \bolds{\beta}, \mathcal{P}) \ge I(P_0 | \bolds
{\beta}, \mathcal{P}_{\mathrm{PL}})$
implies $I(P_0 | \bolds{\beta}, \mathcal{P})^{-1} \le I(P_0 | \bolds
{\beta}, \mathcal{P}_{\mathrm{PL}})^{-1}$, with equality holding if $E_0(X_j|\mathbf{Z}=\mathbf{z})$ are all
additive.

\begin{theorem}\label{thm:eff}
Suppose $I(P_0 | \bolds{\beta}, \mathcal{P}_{\mathrm{PL}})$ is positive
definite. Then,
$I(P_0 | \bolds{\beta}, \mathcal{P})^{-1} < I(P_0 | \bolds{\beta},
\mathcal{P}_{\mathrm{PL}})^{-1}$ unless
$E_0[\bolds{\eta}(\mathbf{Z})-E_0(\mathbf{X}|\mathbf{Z})][\bolds
{\eta}(\mathbf{Z})-E_0(\mathbf{X}|\mathbf{Z})]^\top
=\mathbf{O}$,
where $\mathbf{O}$ is the $p\times p$ matrix with all entries being zero,
and $A<B$ means that $B-A$ is non-negative definite and $A\neq B$.
\end{theorem}

Theorem~\ref{thm:eff} tells that using an additive model for the
nonparametric part can lead to drastic gains of efficiency in the
estimation of the parametric components. The efficiency gains occur if
the parametric covariates $\mathbf{X}$ are approximated by non-additive
transformations of the nonparametric covariates $\mathbf{Z}$. If the
approximation is exact, then estimation of the parametric part in the
partially linear model (\ref{plmodel}) breaks down since $I(P_0 |
\bolds{\beta}, \mathcal{P}_{\mathrm{PL}})=\mathbf{O}$, while it does not
with the partially
linear additive model (\ref{pladdmodel}). If the approximation is very
crude, one has large efficiency gains by using additive models for the
nonparametric part.

\section{Semi-parametric efficient estimation}\label{sec:epa}

Let $\bolds{\beta}^0$ and $m^0$ denote the true parameter values. In this
section we present the semi-parametric efficient estimator of
$\bolds{\beta}^0$ that achieves the minimal asymptotic variance
$I(P_0 |
\bolds{\beta}, \mathcal{P})^{-1}$. The construction is based on a
smooth backfitting
technique and a profiling method. The latter is basically for
estimating the least favorable curve, and is applied to the Gaussian
error model to produce an initial estimator of $\bolds{\beta}^0$ to
be used in
the construction of the semi-parametric efficient estimator.

\subsection{Smooth backfitting methods}\label{subsec:sbf}

The smooth backfitting method, introduced by Mammen, Linton and Nielsen
\cite{mln1999}, is known to be a powerful technique for estimating additive
regression functions. Since our profiling method involves smooth
backfitting for non-additive functions, we discuss some properties of
the method when the target function is not additive.

Let $W$ be a random variable and $\{W^i\}$ be a random sample
distributed as $W$. The smooth backfitting estimator, $\hat m_W^{\mathrm{add}}(\mathbf{z}) \equiv\hat m_{W,0}^{\mathrm{add}}\,+\,\hat m_{W,1}^{\mathrm{add}}(z_1)\,+\,\cdots\,+\,\hat m_{W,d}^{\mathrm{add}}(z_d)$, with responses $W^i$ and
regressors $\mathbf{Z}^i$, are defined as the solution of following integral
equations:
\begin{eqnarray}\label{inteq}
\hat m_{W,j}^{\mathrm{add}} = \tilde{m}_{W,j} - \sum_{l=1, \neq j}^d
\hat\Pi_{j}(\hat{m}_{W,l}^{\mathrm{add}}) - \hat{m}_{W,0}^{\mathrm{add}},
\qquad1
\le j \le d,
\end{eqnarray}
with the constraints $\langle\hat{m}_{W,j}^{\mathrm{add}} , {\mathbf{1}}
\rangle
= 0$ for $1 \le j \le d$. Here, $\hat{m}_{W,0}^{\mathrm{add}}=n^{-1}\sum_{i=1}^n W^i$ and $\tilde m_{W,j}(z_j)$ denotes the
marginal regression kernel estimator obtained by regressing $W^i$ on
$Z_j^i$ only. The operator $\hat\Pi_{j}$ stands for a projection onto a
Hilbert space equipped with a scalar product $\langle\cdot,\cdot
\rangle$; see \cite{ymp2008} for details. For example, in
the case where $\tilde m_{W,j}(z_j)$ are the local constant marginal
estimators, $\langle g,h \rangle=\int g(\mathbf{z})h(\mathbf{z}) \hat
q_{\mathbf{Z}
}(\mathbf{z})\,\mathrm{d}\mathbf{z},$ with $ \hat q_{\mathbf{Z}}(\cdot)$ being the kernel
estimator of the
design density $q_{\mathbf{Z}}$. Smoothing to the direction of $Z_j$
is done
by the boundary corrected kernel $K_{h_j}(u,v) = c_j(v)h_j^{-1}K^0
((u-v)/h_j)$, where $K^0$ is a base kernel function, $h_j$ is the
bandwidth, and $c_j(v)$ is a factor that gives $\int K_{h_j}(u,v)\,\mathrm{d}u=1$.

Let $m_W(\mathbf{z})=E(W|\mathbf{Z}=\mathbf{z})$. We do not assume that
$m_W$ is an
additive function. Define $m_W^{\mathrm{add}}=m_{W,1}^{\mathrm{add}}+\cdots
+m_{W,d}^{\mathrm{add}}$ to be the projection of $m_W$ onto the space of
additive functions $\mathcal{H}(q_{\mathbf{Z}})$. Then,
$E[m_W(\mathbf{Z})-E(W)-m_W^{\mathrm{add}}(\mathbf{Z})]\delta(\mathbf{Z})=0$ for any $\delta\in\mathcal
{H}(q_{\mathbf{Z}})$. The
additive function $m_W^{\mathrm{add}}(\mathbf{z})$ plays the role of the target
function that the smooth backfitting estimator $\hat m_W^{\mathrm{add}}(\mathbf{z})$ aims at. Lu \textit{et al.} \cite{llty2007} discussed the property of the
smooth backfitting estimators under non-additive regression models in
the context of spatial data analysis. However, they treated only the
case where the bandwidth is asymptotic to $n^{-1/5}$. Below, we give a
uniform expansion of the smooth backfitting estimator for a wider range
of the bandwidths, after tedious asymptotic calculation following the
lines of the arguments in \cite{mln1999}. To state the
theorem, let $\varepsilon=W-E(W)-m_W^{\mathrm{add}}(\mathbf{Z})$ and define
$\varepsilon^i$ accordingly. Let $\tilde m_{\varepsilon,j}(z_j)$ and
$\tilde m_{\varepsilon,j}^{LL}(z_j)$ denote, respectively, the local
constant and linear estimators with responses $\varepsilon^i$ and the
scalar regressors $Z_j^i$. Let $h_j$ be the bandwidth associated with
$Z_j$. The theorem relies on the following assumptions.

\renewcommand{\theassu}{A}
\begin{assu}
\begin{enumerate}[(A1.)]
\item[A1.] For $1 \le j\neq k \le d$, $q_{Z_j,Z_k}$ are bounded away
from zero and infinity on its
support, $[0,1]^2$, and have continuous partial derivatives.
\item[A2.] The base kernel function $K^0$ is symmetric, supported on a
compact support and has bounded derivative.
\item[A3.] The functions $m_{W,j}^{\mathrm{add}}$'s are twice continuously
differentiable.
\item[A4.] $E|W-m_W(\mathbf{Z})|^{r_0}<\infty$ for some $r_0>5/2$.
\end{enumerate}
\end{assu}

\begin{theorem}\label{prop:unifsbf}
Assume that the conditions \textup{A}1--\textup{A}4 hold, and that $h_j$ are asymptotic
to $n^{-\alpha}$ for $1/5 \leq\alpha< 1/2$. Then, for $1 \le j \le
d,$ it holds that
\[
\sup_{z_j \in[0,1]} |\hat m_{W,j}^{\mathrm{add}}(z_j)-m_{W,j}^{\mathrm{add}}(z_j)-h_j a_{1,j,n}(z_j)-h_j^2a_{2,j}(z_j)-\tilde
m_{\varepsilon,j}(z_j)|=\mathrm{o}_p((nh_j)^{-1/2})
\]
in the local constant case, and that
\[
\sup_{z_j \in[0,1]} |\hat m_{W,j}^{\mathrm{add}}(z_j)-m_{W,j}^{\mathrm{add}}(z_j)-h_j^2a_{3,j}(z_j)-\tilde
m_{\varepsilon,j}^{LL}(z_j)|=\mathrm{o}_p((nh_j)^{-1/2})
\]
in the local linear case, for some functions $a_{1,j,n}$ that are
uniformly bounded and non-zero only for $z_j \in[0,ch_j) \cup
(1-ch_j,1]$ for some constant $0<c<\infty$, and for some functions
$a_{2,j}$ and $a_{3,j}$ that are continuous.
\end{theorem}

A proof of Theorem~\ref{prop:unifsbf} can be found in an extended
version of this paper that can be downloaded from \url{http://stat.snu.ac.kr/theostat/papers/BEJ296\_ExtendedVersion.pdf}.

\subsection{Profiling with Gaussian error models}\label{subsec:profile}

We apply a profiling technique to remove the infinite-dimensional
parameter $m$ in the estimation of $\bolds{\beta}^0$. For a general framework
of profiling approaches to semi-parametric models, we refer to \cite{sw1992}. See also \cite{madv2000} for a more
recent work on profile likelihood.

Define $\hat m_{\mathbf{X}}^{\mathrm{add}}=(\hat m_{X_1}^{\mathrm{add}}, \ldots
, \hat
m_{X_p}^{\mathrm{add}})^\top$. We note that $\hat m_{\mathbf{X}}^{\mathrm{add}}$ is an
estimator of $\bolds{\eta}$ and $\hat m_{Y}^{\mathrm{add}}$ is an
estimator of
$\bolds{\beta}^{0\top}\bolds{\eta}+ m^0$. For each given $\bolds
{\beta}$, let $\hat
m^{\mathrm{add}}(\mathbf{z};\bolds{\beta})=\sum_{j=1}^d \hat
m_j^{\mathrm{add}}(z_j;\bolds{\beta})$\vspace{1pt}
be the
smooth backfitting estimator obtained by taking
$Y^i-\mathbf{X}^{i\top}\bolds{\beta}=\mathbf{X}^{i\top}(\bolds
{\beta}^0-\bolds{\beta})+m^0(\mathbf{Z}
^i)+\epsilon^i$
as responses and $\mathbf{Z}^i$ as covariates. Recall that the least favorable
curve is given by $m^*(\cdot, \bolds{\beta}) \equiv\bolds{\eta
}^\top
(\bolds{\beta}^0-\bolds{\beta}) + m^0$. Thus, we may regard $\hat
m^{\mathrm{add}}(\cdot;\bolds{\beta})$\vspace{2pt} as an estimator of the least favorable curve
$m^*(\cdot, \bolds{\beta})$. Since $\hat m^{\mathrm{add}}(\mathbf{z};\bolds
{\beta})=\hat
m_{Y}^{\mathrm{add}}(\mathbf{z})-\hat m_{\mathbf{X}}^{\mathrm{add}}(\mathbf{z})^\top
\bolds{\beta}$\vspace{2pt}
by the fact that the smooth backfitting operation is linear in response
vectors, the estimated profile likelihood based on the Gaussian error
model is given by
\[
-\sum_{i=1}^n [Y^i-\mathbf{X}^{i\top} \bolds{\beta}- \hat m^{\mathrm{add}}(\mathbf{Z}
^i;\bolds{\beta})
]^2 = -\sum_{i=1}^n \bigl[ Y^i-\hat m_{Y}^{\mathrm{add}}(\mathbf{Z}^i)-\bigl(\mathbf
{X}^i -\hat
m_{\mathbf{X}}^{\mathrm{add}}(\mathbf{Z}^i) \bigr)^\top\bolds{\beta}\bigr]^2.
\]
The estimator that maximizes the above Gaussian profile likelihood is
then given by
\[
\hat{\bolds{\beta}}= \Biggl( \sum_{i=1}^n \tilde\mathbf{X}^i \tilde
\mathbf{X}^{i\top} \Biggr)^{-1} \Biggl(
\sum_{i=1}^n\tilde\mathbf{X}^i \tilde Y^i \Biggr),
\]
where $\tilde\mathbf{X}^i= \mathbf{X}^i -\hat m_{\mathbf{X}}^{\mathrm{add}}(\mathbf{Z}^i)$ and
$\tilde
Y^i= Y^i - \hat m_{Y}^{\mathrm{add}}(\mathbf{Z}^i)$.

\begin{theorem}\label{thm:asymprofile}
Suppose that the assumptions \textup{A}1--\textup{A}4 hold with $W=Y$ and $X_j$, $1 \le j
\le p$. Also, assume that $E[\exp(|X_j-E(X_j|\mathbf{Z})|)|\mathbf
{Z}]<C$ a.s. for
some $C>0$, $1 \le j \le p$. If the bandwidths $h_j$ are asymptotic to
$n^{-\alpha}$ for $1/5 \leq\alpha<1/2$, then it holds that
\[
\sqrt n ( \hat{\bolds{\beta}}-\bolds{\beta}^0) \stackrel d
\Rightarrow N\bigl(\mathbf{0},\operatorname{var}(\epsilon) \bigl[E \bigl(\mathbf{X}-\bolds{\eta}(\mathbf{Z})\bigr)\bigl(\mathbf
{X}-\bolds{\eta}(\mathbf{Z})\bigr)^\top\bigr]^{-1}\bigr).
\]
\end{theorem}

A proof of Theorem~\ref{thm:asymprofile} is given in the \hyperref[app]{Appendix}. We
note that the asymptotic variance of the estimator $\hat{\bolds{\beta}}$ is
larger than $I(P_0 | \bolds{\beta}, \mathcal{P})^{-1}$. This can be
seen directly from
a projection property. In fact, $\operatorname{var}(\epsilon) \ge I_g^{-1}$ and
the equality hold if $g$ is Gaussian. This means that the estimator
$\hat{\bolds{\beta}}$ achieves the semi-parametric efficiency in the reduced
model where $g$ is specified as a Gaussian density. It is also
interesting to see what happens if $\eta_0(\mathbf{X},\mathbf{Z})
\equiv
E_0(Y|\mathbf{X},\mathbf{Z})$ does not belong to the partially linear
additive model
of the form (\ref{pladdmodel}). In this case, our estimator of $\eta_0$
converges to $\eta^*$, which is the $L_2(q)$-projection of $\eta_0$
onto the space
\begin{equation}\label{modelspace}
\mathcal{F}= \{ f \in L_2(q) \mid f(\mathbf{x},\mathbf{z}) = \bolds{\beta
}^{\top} \mathbf{x}+
m(\mathbf{z}),\ \bolds{\beta}
\in\mathbb{R}^p,\ m \in\mathcal{H}\}.
\end{equation}

\subsection{Adapting to unknown error density}\label{subsec:adaptive}

In this subsection, we construct the semi-parametric efficient
estimator that achieves the minimal asymptotic variance discussed in
Section \ref{sec:eff}. We follow the approach adopted by Bickel \cite{b1982},
Schick \cite{s1986,s1993}, Park \cite{p1990}, Cuzick \cite{c1992} and Bhattacharya and
Zhao \cite{bz1997}. Write $I= I(P_0|\bolds{\beta},\mathcal{P})$ and define
$\bolds{\beta}_n^* =
\bolds{\beta}^0 - I^{-1} n^{-1} \sum_{i=1}^n [\mathbf{X}^i- \eta
(\mathbf{Z}^i)]
\varphi(\epsilon)$. Then, the random sequence $\bolds{\beta}_n^*$
achieves the
efficiency bound. We plug some estimators of the unknown quantities
into $\bolds{\beta}_n^*$. We estimate the error density $g$ by using the
`pseudo' errors $\hat\epsilon^i \equiv\tilde Y^i - \tilde
\mathbf{X}^{i\top}\hat{\bolds{\beta}}$, where $\hat{\bolds{\beta}}$
is the Gaussian profile
estimator constructed in Section~\ref{subsec:profile}. In particular,
we take $\hat g(t)=b+(na)^{-1}\sum_{i=1}^n L((t- \hat\epsilon^i)/a)$
and $\hat g'(t)=\mathrm{d}\hat g(t)/\mathrm{d}t$, where $a$ and $b$ are positive
constants that depend on the sample size $n$, and $L$ is a symmetric
differentiable density function. Define
\[
\hat I=\Biggl(n^{-1} \sum_{i=1}^n \tilde\mathbf{X}^i\tilde\mathbf
{X}^{i\top}\Biggr)\Biggl(
n^{-1}\sum_{i=1}^n \hat\varphi(\hat\epsilon^i)^2\Biggr),
\]
where $\hat\varphi$ is the `symmetrized' estimator of $\varphi$
defined by $\hat\varphi(e) = [ (\hat g'/\hat g)(e) - (\hat g'/\hat
g)(-e)]/2$. Our semi-parametric efficient estimator is then given by
\[
\tilde{\bolds{\beta}}= \hat{\bolds{\beta}}- \hat I^{-1} \frac1 n
\sum_{i=1}^n
\tilde\mathbf{X}^i \hat\varphi(\hat\epsilon^i).
\]

\renewcommand{\theassu}{B}
\begin{assu}
\begin{enumerate}[(B1.)]
\item[B1.] The error $\epsilon$ has an absolutely continuous and
symmetric density $g$ with
respect to the Lebesgue measure, $\mu$, and $I_{g} = \int(g'^2 /g)\,\mathrm{d}\mu<\infty$.
\item[B2.] The kernel $L$ is a symmetric density function with three
bounded and Lipschitz continuous derivatives.
\item[B3.] The sequences $a$ and $b$ converge to zero, as $n
\rightarrow\infty$, and satisfy $n^{1/2}h_j b(a^2 \wedge b^2)
\rightarrow\infty$ and $a^2 / \{h_j (\log n )^2\} \rightarrow\infty$
for all $1 \le j \le d$.
\end{enumerate}
\end{assu}

\begin{theorem}\label{thm:asymadap}
Assume that the conditions of Theorem \ref{thm:asymprofile} and the
assumptions \textup{B}1--\textup{B}3 hold. Then, $\sqrt n ( \tilde{\bolds{\beta}}-\bolds
{\beta}^0)
\stackrel d \Rightarrow N({\mathbf{0}},I(P_0 | \bolds{\beta}, \mathcal{P})^{-1})$.
\end{theorem}

A proof of Theorem~\ref{thm:asymadap} is given in the \hyperref[app]{Appendix}. For a
choice of the bandwidth $a$ in $\hat{g}$, one can devise a data-driven
choice along the lines of Park \cite{p1993}. For $h$, one can follow the
approach of Mammen and Park \cite{mp2005}. In this adaptation step,
misspecification of the model may result in a meaningless estimator.
This is in contrast to the estimation in the initial step where the
procedure estimates the projection of the mean function onto the model
space $\mathcal{F}$ at (\ref{modelspace}). The reason is that the residuals
from the initial step include not only the pure errors but also the
deviation of the true regression function from its projection onto
$\mathcal{F}$. These residuals mislead estimation of the score function.

\section{Numerical properties}\label{sec:sim}

We generated 500 random samples of the size $n=400$. We used
Epanechnikov kernel for the regression and the Gaussian density kernel
for the estimation of the score function. We applied a local constant
version of smooth backfitting. We took $m_1(z_1)
=\sin\{2\uppi(z_1-0.5)\}$ and $m_2(z_2) =z_2-0.5+\sin\{2\uppi(z_2-0.5)\}$.
We set $m_0=3$, $\beta_1=1.5$ and $\beta_2=0.8$. We drew $(Z_1,Z_2)$
from $N_2((0.5,0.5)^\top, \Sigma)$ truncated to $[0,1]^2$, where
$\Sigma= \{(1-\rho)I + \rho{\mathbf{1}\mathbf{1}}^\top\}/4 $. We generated
$X_1=C Z_1(1-2Z_2)+U$ for some constant $C$, where $U \sim N(0, 0.5)$,
and $X_2$ from $\mbox{Bernoulli}(p(X_1, Z_1,Z_2))$, where $p(X_1,
Z_1,Z_2)= g(\exp((Z_1+Z_2)/2)+\sin(2\uppi Z_1)-X_1^2)$ and $g(t)=
\exp(t)/(1+\exp(t))$. Note that $E(X_1|\mathbf{Z}=\cdot)$ is
orthogonal to the
space of additive functions.

We compared the Gaussian profile estimator (SAM), given in Section
\ref{subsec:profile}, and the profile kernel estimator (PL), given in
\cite{s1988}, which is for the partial linear model without the
additive structure. For this, we generated $\epsilon$ from $N(0,1)$ and
set $\rho=0$. In the case where $p=1$, that is, $X_2$ does not enter
the model, the theoretical value of the ratio of the asymptotic
variance of SAM to that of PL equals $1/(1+0.1707C^2)$. The empirical
values from our simulation study for the bandwidth pair $(h_1, h_2)$
that gave the best mean square error (MSE) were $0.7818, 0.5868$ and
$0.4082$ for $C=1,2$ and $3$, respectively, which nearly coincided with
the theoretical values. We tried other values of $\rho$, but the lesson
was the same. In the case where $p=2$ and $d=5$ with $(Z_1,\ldots
,Z_5)$ from $N_5((0.5,\ldots,0.5)^\top, \Sigma)$ truncated to
$[0,1]^5$ and $m_j(z_j)=z_j^2$ for $3 \le j \le5$, we took $C=1$ and
found that SAM beat PL for all bandwidth choices that we tried. The
Gaussian profile estimator was stable while PL broke down for small
bandwidths. The best MSE of SAM and that of PL, respectively, for
various choices of the bandwidth pair $(h_1, h_2)$ were $0.0032$ and
$0.0051$ for $\beta_1$ and $0.0186$ and $0.0269$ for $\beta_2$.

\begin{figure}

\includegraphics{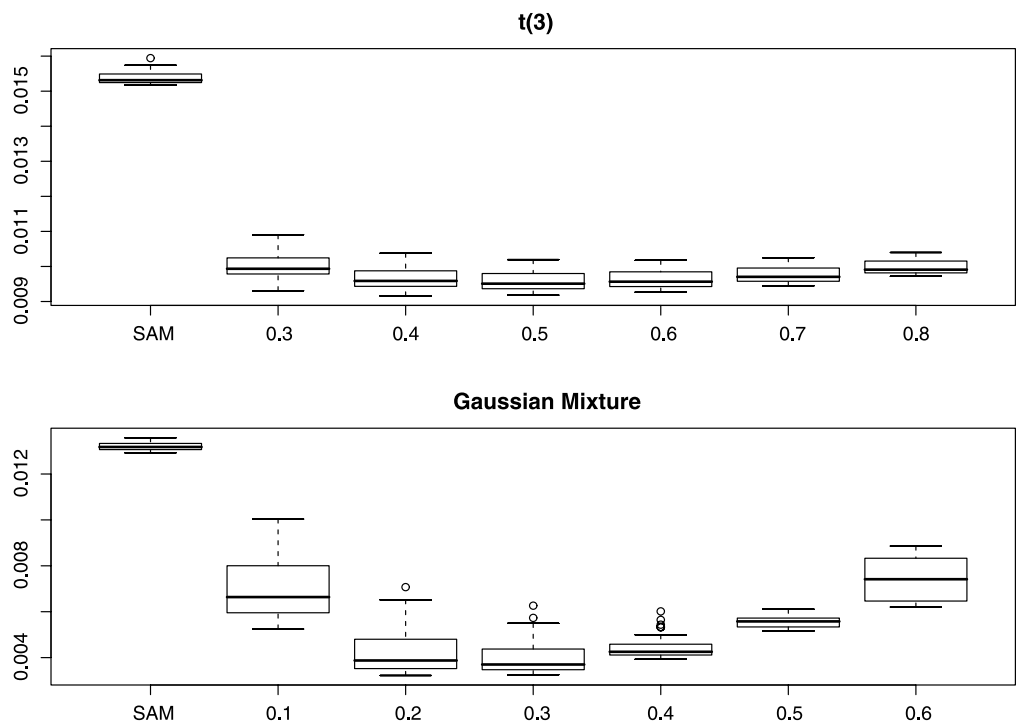}

\caption{Mean square errors of SAM and ASAM.}
\label{fig:eff1}\vspace*{-15pt}
\end{figure}

Next, we compared SAM with the semi-parametric efficient estimator
(ASAM). For this, we considered the case where $p=d=2, C=1$ and
$\rho=0.8$, and generated $\epsilon$ from $N(0,1)$, $t$-distribution
with degree of freedom 3, and $\frac{1}{2}N(-1.5, 0.6^2) +
\frac{1}{2}N(1.5, 0.6^2)$. For ASAM, we took $b=0.01$, and six
different choices of $a$: $a_i =0.3+0.1i, 0 \le i \le5$, for $N(0,1)$
and $t(3)$ errors and $a_i =0.1+0.1i, 0 \le i \le5$, for the Gaussian
mixture error. We used $36$ different choices for the bandwidth pair
$(h_1,h_2) \in\{0.05,0.10,,\ldots, 0.30\}^2$. Figure~\ref{fig:eff1} is
for the estimators of $\beta_1$. Each box-plot was obtained from the
$36$ values of MSE that corresponded to the $36$ bandwidth pairs
$(h_1,h_2)$. For ASAM, the value of $a$ is indicated on the horizontal
scale. The figure suggests that the values of the MSE of ASAM are far
smaller than those of SAM for the entire range of the bandwidth $a$,
under $t(3)$ and the Gaussian mixture error models. The box-plots for
the Gaussian error model are not given here since SAM and ASAM gave
similar performance. The results for $\beta_2$ are not reported either
since they give a similar lesson.\vspace*{-2pt}

\section{Boston housing data}\label{sec:boston}\vspace*{-1pt}

We applied the semi-parametric efficient estimators to Boston housing
data as an illustration. As in \cite{wy2007,fh2005}, we took the median price in 1,000 USD (MEDV) as the response
$Y$. Also, we chose\vadjust{\goodbreak} as covariates $X_1$, $X_2$ and $Z_1, \ldots, Z_6$,
respectively, the eight variables LSTAT (percentage values of lower
status population), CHAS (a dummy variable that takes the value $1$ if
the tract borders Charles River; 0 otherwise), CRIM (per capita crime
rate), RM (average numbers of rooms per dwelling), NOX (nitric oxides
concentration), PTRATIO (pupil--teacher ratios), DIS (weighted distances
to five Boston employment centers) and TAX (full-value property tax
rate per 10,000~USD). The logarithms of LSTAT, DIS and TAX were taken
to reduce sparse areas, as in \cite{wy2007}. We chose the model
$Y = m_0 + \beta_1 X_1 +\beta_2 X_2 +\sum_{j=1}^6 m_j (Z_j)+\epsilon$.
In the data set, there were 16 cases for which $Y$ took the maximal
value $50$. These may be censored responses that one may remove from
analysis. Indeed, an initial analysis showed a strong asymmetry in the
distribution of the residuals, which led us to exclude the 16 cases for
further analysis. For additive regression, we applied local constant
smooth backfitting with the Epanechnikov kernel and bandwidths $h_j$
chosen by a rule of thumb.

With SAM, we obtained $\hat\beta_1= -6.203$ and $\hat\beta_2= 0.985$.
Their estimated standard errors were $0.420$ and $0.597$, respectively.
This suggests that $\hat\beta_2$ is not strongly significant while
$\hat\beta_1$ is. The generalized $R^2$ was $0.862$. For ASAM, in the
estimation of the score function, we used a bandwidth $a$ that was
obtained by R function \texttt{bw.SJ()}. With ASAM, we got $\tilde
\beta_1= -6.172$ and $\tilde\beta_2= 1.366$, and their estimated
standard errors were $0.399$ and $0.567$, respectively. Thus, with
ASAM, both the estimated coefficients are strongly significant. This
may be an indication that a Gaussian error model is not appropriate for
the data set. The generalized $R^2$ was almost the same as in the
analysis with SAM.

\begin{appendix}\label{app}
\section*{Appendix}

\begin{pf*}{Proof of Theorem \ref{thm:asymprofile}}
We only treat
the case with local constant smooth backfitting. The case with local
linear smooth backfitting can be dealt with similarly. We prove
\renewcommand{\theequation}{\arabic{equation}}
\begin{equation}\label{pfconv1}
n^{-1/2} \sum_{i=1}^n \tilde\mathbf{X}^i ( \tilde Y^i -\tilde
\mathbf{X}^{i\top}\bolds{\beta}^0 ) - n^{-1/2} \sum_{i=1}^n
\bigl(\mathbf{X}^i
-\bolds{\eta}(\mathbf{Z}^i)\bigr)\epsilon^i = \mathrm{o}_p(1).
\end{equation}
Write $\Delta(\mathbf{z}) = m^0(\mathbf{z})-\hat m^{\mathrm{add}}(\mathbf{z};\bolds
{\beta}
^0)$. The
left-hand side of equation (\ref{pfconv1}) equals $C_1 + C_2 + C_3$,
where $C_1=n^{-1/2} \sum_{i=1}^n (\mathbf{X}^i -\bolds{\eta
}(\mathbf{Z}^i))\Delta(\mathbf{Z}^i)$,
$C_2=n^{-1/2} \sum_{i=1}^n (\bolds{\eta}(\mathbf{Z}^i)-\hat
m_{\mathbf{X}}^{\mathrm{add}}(\mathbf{Z}^i))\epsilon^i$, and $C_3=n^{-1/2} \sum_{i=1}^n
(\bolds{\eta}(\mathbf{Z}^i)-\hat m_{\mathbf{X}}^{\mathrm{add}}(\mathbf
{Z}^i))\Delta(\mathbf{Z}^i)$. Write
$\Delta(\mathbf{z}) = \Delta_0 + \sum_{j=1}^d \Delta_j(z_j)$. By
Theorem~\ref{prop:unifsbf}, standard techniques of kernel smoothing,
integration by part and the representation of $m^0$ and $\hat m^{\mathrm{add}}(\mathbf{z};\bolds{\beta}^0)$ as a solution of an integral equation with
differentiable kernel (see equation (\ref{inteq})), we have
\[
\sup_{\mathbf{z}\in[0,1]^d}|\Delta(\mathbf{z})|=\mathrm{o}_p(\delta_n),
\qquad
\sup_{z_j \in[0,1]} \biggl|\frac {\mathrm{d}} {\mathrm{d}z_j }\Delta_j(z_j)- h_j
b_{n,j}(z_j)\biggr|=\mathrm{o}_p(\delta_n)
\]
for some uniformly bounded non-random functions $b_{n,j}$, where
$\delta_n=n^{-a}$ for some $a\in(0, 1/2-\alpha)$. These imply that
$\delta_n^{-1}\Delta\in B(\mathbf{0},1)$ with probability tending to one,
where $B(\mathbf{0},1)$ denotes a class of additive functions
$\sum_{j=1}^{d} g_j(z_j)$ such that each $g_j$ is a real function
defined on $[0,1]$ and satisfies $\sup_{t, t' \in[0,1]} |g_j(t)-
g_j(t')|\leq|t-t'|$. The covering number with bracketing of $B(\mathbf{0},1)$ with respect to sup-norm, $N_{[\cdot]}(\eta) \equiv
N_{[\cdot]}(\eta,B(\mathbf{0},1),\| \cdot\|_{\infty})$, is bounded by
$(2\eta^{-1})^{d}3^{d \eta^{-1}}$. Define random functionals
$F(X_j^i,\mathbf{Z}^i)\dvtx B(\mathbf{0},1) \rightarrow\mathbb{R}$ by
$[F(X_j^i,\mathbf{Z}^i)](g)=(X_j^i -\eta_j(\mathbf{Z}^i))g(\mathbf
{Z}^i)$, and $F_j\dvtx
B(\mathbf{0},1) \rightarrow\mathbb{R}$ by $F_j = n^{-1/2}\sum_{i=1}^n
F(X_j^i,\mathbf{Z}^i)$. Then, using Corollary 8.8 of van de Geer
\cite{vdg2000} and
the tail condition assumed in the theorem, one can show $\sup_{g \in
B(\mathbf{0},1)} |F_j g|=\mathrm{O}_p(1)$. Let $C_{1,j}$ denote the $j$th element
of $C_1$. Since $P(|\delta_n^{-1}C_{1,j}| > M )\le P(\sup_{g \in
B(\mathbf{0},1)} |F_j g| > M) + P(\delta_n^{-1}\Delta\notin B(\mathbf{0},1))$, we
obtain $C_{1,j}=\mathrm{O}_p(\delta_n)=\mathrm{o}_p(1)$. One can prove $C_{2}=\mathrm{o}_p(1)$
using a truncation argument with Theorem~\ref{prop:unifsbf} and
applying the Chebyshev inequality conditioning on $(\mathbf{X}^i,
\mathbf{Z}^i)$. The
fact that $C_{3}=\mathrm{o}_p(1)$ follows from $P(Z_j^i$ lies in $[0,ch_j)\cup(1-ch_j, 1])=\mathrm{O}(h_j)$ for some constant $0<c<\infty$ and
Theorem~\ref{prop:unifsbf}.
\end{pf*}

\begin{pf*}{Proof of Theorem \ref{thm:asymadap}}
We will show that
$\tilde{\bolds{\beta}}- \bolds{\beta}_n^* = \mathrm{o}_p(n^{-1/2})$. It
suffices to show
\renewcommand{\theequation}{\arabic{equation}}
\begin{equation}\label{mclaim}
\hat I^{-1} n^{-1} \sum_{i=1}^n \tilde\mathbf{X}^i \hat\varphi
(\hat
\epsilon^i) = \hat{\bolds{\beta}}-\bolds{\beta}^0 +I^{-1} n^{-1}
\sum_{i=1}^n [\mathbf{X}^i-
\eta(\mathbf{Z}^i)] \varphi(\epsilon^i)+\mathrm{o}_p(n^{-1/2}).
\end{equation}
By Theorem~\ref{prop:unifsbf} and standard techniques of kernel
smoothing along with assumption B3, it holds that, uniformly over $i$,
\begin{equation}\label{exp_phi}
\hat\varphi(\hat\epsilon^i)=\hat\varphi( \epsilon^i)- \tilde
\mathbf{X}^{i\top} (\hat{\bolds{\beta}}-\bolds{\beta}^0)\hat
\varphi'(\epsilon^i)- \{
\hat
m^{\mathrm{add}} (\mathbf{Z}^i;\hat{\bolds{\beta}})-m^0(\mathbf{Z}^i)\}
\hat
\varphi'(\epsilon^i)+\mathrm{o}_p(n^{-1/2}).
\end{equation}
Also, using the proof of Lemma~4.1 in \cite{b1982} and standard
calculus, one can show $\hat I = I+o_p(1)$ and $n^{-1}\sum_{i=1}^n
\tilde\mathbf{X}^i \tilde\mathbf{X}^{i\top} \hat\varphi
'(\epsilon^i)= -I+\mathrm{o}_p(1)$.
Thus, the proof of the theorem is completed if we verify
\begin{eqnarray}
\label{effclaim3}n^{-1} \sum_{i=1}^n \tilde\mathbf{X}^i \{ \hat m^{\mathrm{add}} (\mathbf
{Z}^i;\hat
{\bolds{\beta}})-m^0(\mathbf{Z}^i)\}\hat\varphi'(\epsilon^i) &=&
\mathrm{o}_p(n^{-1/2});
 \\
\label{effclaim4}n^{-1} \sum_{i=1}^n \tilde\mathbf{X}^i \hat\varphi(\epsilon^i) -n^{-1}
\sum_{i=1}^n \{\mathbf{X}^i- \bolds{\eta}(\mathbf{Z}^i)\} \varphi
(\epsilon^i)& =&
\mathrm{o}_p(n^{-1/2}).
\end{eqnarray}
Proofs of (\ref{effclaim3}) and (\ref{effclaim4}) can be based on the
following lemma, which follows from Corollary 2.7.4 in \cite{vdvw1996} and assumption B2 on $L$. Note that the moment
condition on $\epsilon$ ensures the entropy bound. To state the lemma,
define
\[
\mathcal{C}_M^{\alpha}(\mathcal{X})=\biggl\{f\dvtx\mathcal{X} \rightarrow
\mathbb{R}\dvtx \sup_x|f(x)| + \sup_{x,y}\frac{ |f(x)-f(y)|^{\alpha}
}{|x-y|}\leq M \biggr\}
\]
for a set $\mathcal{X}\subset\mathbb{R}$ and a real number $\alpha
\in
(0,1]$. Let $\| \cdot\|_{g}$ denote the $L_2$ norm with respect to the
density $g$.
\begin{lemma}\label{lemma:entropy}
Assume the conditions of Theorem \ref{thm:asymadap}. Then there exists
a constant $M$ such that, with probability tending to one, $b(a \wedge
b)\hat\varphi\in\mathcal{C}_{M}^{1}(\mathbb{R})$, $[nh_{\max}a^6b/ (\log n)^2]^{1/2}(\hat\varphi-\varphi_n) \in
\mathcal{C}_{M}^{1}(\mathbb{R})$ and $b(a^2 \wedge b^2)\hat\varphi'
\in\mathcal{C}_{M}^{1}(\mathbb{R})$. Moreover, there exist constants
$\delta>0$ and $C_1>0$ such that $\log
N_{[\cdot]}(\eta,\mathcal{C}_{M}^{1}(\mathbb{R}),\| \cdot\|_{g})
\leq
C_1 \eta^{-(2-\delta)}$.
\end{lemma}
\upqed
\end{pf*}
\end{appendix}

\section*{Acknowledgement}
Research of Kyusang Yu was supported in part by Basic Science Research Program
through the National Research Foundation of Korea (NRF) funded by the Ministry of Education,
Science and Technology (2010-0023488). Research of Byeong U. Park was supported by the Mid-career Researcher
Program through NRF grant funded by the MEST (No. 2010-0017437).

\printhistory

\end{document}